\documentstyle[12pt,amsfonts]{article}
\def\reducedqp{{\mathbb C}_{rq}}
\def\reducedqpk{{\mathbb C}^{k}_{rq}}

\def\reducedqpzero{{\mathbb C}^{0}_{rq}}
\def\reducedqpone{{\mathbb C}^{1}_{rq}}
\begin{document}
\title{On a graded $q$-differential algebra}
\author{Viktor Abramov}%
\date{}%
\maketitle%
\begin{abstract}%
We construct the graded $q$-differential algebra on a ${\mathbb
Z}_N$-graded algebra by means of a graded $q$-commutator. We apply
this construction to a reduced quantum plane and study the first
order differential calculus on a reduced quantum plane induced by
the $N$-differential of the graded $q$-differential algebra.
\end{abstract}%
\section{Graded $q$-differential algebra}%
In this section given a ${\mathbb Z}_N$-graded algebra we
construct the graded $q$-dif\-fer\-ential algebra. Let us remind
the definition of a graded $q$-differential algebra
(\cite{Dubois1}). A unital associative algebra is said to be a
graded $q$-differential algebra ($q\in {\mathbb C}, q\neq 1$) if
it is a ${\mathbb Z}_N$-graded (or $\mathbb Z$-graded) algebra
endowed with the linear mapping $d$ of degree $+1$ satisfying the
graded $q$-Leibniz rule and $d^N=0$ in the case when $q$ is a
primitive $N$-th root of unity. The linear mapping $d$ is called
an $N$-differential of a graded $q$-differential algebra.

 Let $\cal A$ be an associative unital
$\mathbb Z$ (or ${\mathbb Z}_N$)-graded algebra over the complex
numbers $\mathbb C$ and ${\cal A}^k\subset {\cal A}$ be the
subspace of homogeneous elements of a grading $k$. The grading of
a homogeneous element $w$ will be denoted by $|w|$, which means
that if $w\in {\cal A}^k$ then $|w|=k$. Let $q$ be a complex
number such that $q\neq 1$. The $q$-commutator of two homogeneous
elements $w,w'\in {\cal A}$ is defined by the formula
\begin{equation}
[w,w']_q=ww'-q^{|w||w'|}w'w.
\end{equation}
Using the associativity of an algebra $\cal A$ and the property
$|ww'|=|w|+|w'|$ of its graded structure it is easy to show that
for any homogeneous elements $w,w',w''\in {\cal A}$ the
$q$-commutator has the property
\begin{equation}
[w,w'w'']_q=[w,w']_q w''+q^{|w||w'|}w'[w,w'']_q.%
\label{q-commutator property}
\end{equation}
Given an element $v\in {\cal A}^1$ one can define the mapping
$d_v:{\cal A}^k\to {\cal A}^{k+1}$ by the following formula
$$
d_vw=[v,w]_q,\qquad w\in {\cal A}^k.
$$
It follows from the property of $q$-commutator (\ref{q-commutator
property}) that $d_v$ is a graded $q$-differential on an algebra
$\cal A$, i.e. it is a homogeneous linear mapping of degree 1
satisfying the graded $q$-Leibniz rule%
\begin{equation}
d_v(ww')=d_v(w)w'+q^{|w|}wd_v(w'), %
\label{q-Leibniz}
\end{equation}
where $w,w'$ are the homogeneous elements of ${\cal A}$.

{\bf Lemma 1.1} {\it For any integer $k\geq 2$ the
$k$-th power of
the $q$-differential $d_v$ can be written as follows%
\begin{equation}
d_v^kw=\sum_{i=0}^{k}p_i^{(k)}v^{k-i}wv^i,
\end{equation}
where $w$ is a homogeneous element of $\cal A$ and}%
\begin{eqnarray}
p^{(k)}_i&=&(-1)^iq^{|w|_i}\frac{[k]_q!}{[i]_q![k-i]_q!}=
(-1)^iq^{|w|_i}\left[%
\begin{array}{c}
  k \\
  i \\
\end{array}%
\right]_q,\;\\ |w|_i&=&i|w|+\frac{i(i-1)}{2}.
\end{eqnarray}

This lemma can be proved by means of the mathematical induction
and the following identities
\begin{eqnarray}
p^{(k)}_0&=&p^{(k+1)}_0=1,\;\;\;\;\; p^{(k+1)}_{k+1}=-q^{|w|+k}p^{(k)}_k,\nonumber\\
  p^{(k+1)}_i&=&p^{(k)}_i-q^{|w|+k}p^{(k)}_{i-1},\;\;\; 1\leq i\leq
  k.\nonumber
\end{eqnarray}

{\bf Theorem 1.1} {\it If $N$ is an integer such that $N\geq 2$,
$\cal A$ is a ${\mathbb Z}_N$-graded algebra, $q$ is a primitive
$N$-th root of unity and $v^N=\alpha e$, where $\alpha \in
{\mathbb C}$ and $e$ is the unity element of an algebra $\cal A$,
then $d^N_vw=0$ for any $w\in {\cal A}$.}

It follows from the Lemma 1.1 that if $q$ is a primitive $N$-th
root of unity then for any integer $l=1,2,\ldots,N-1$ the
coefficient $p^{(N)}_l$ contains the factor $[N]_q$ which is equal
to zero in the case of $q$ being a primitive $N$-th root of unity
and this implies $p^{(N)}_l=0$. Thus
$d^N_v(w)=v^Nw+(-1)^Nq^{|w|_N}wv^N$. Taking into account that
$v^N=\alpha e$ we obtain $d^N_v(w)=(1+(-1)^Nq^{|w|_N})\alpha w$.
The first factor in the right-hand side of the above formula
equals to zero. Indeed if $N$ is an odd number then
$1-(q^{N})^{\frac{N-1}{2}}=0$. In the case of an even integer $N$
we have $1+(q^{\frac{N}{2}})^{N-1}=1+(-1)^{N-1}=0$, and this ends
the proof.

Let $\cal A$ be an associative unital ${\mathbb Z}_N$-graded
algebra over the complex numbers $\mathbb C$ with unit element
denoted by $e$. Then from the property (\ref{q-Leibniz}) and the
Theorem 1.1 it follows

{\bf Corollary.} {\it If there exists an element $v\in {\cal A}$
of grading 1 such that $v^N=\alpha e, \alpha\in {\mathbb C}$ then
an algebra $\cal A$ endowed with the homogeneous linear mapping
$d_v:{\cal A}^k\to {\cal A}^{k+1}$ of degree +1, defined by
$d_vw=[v,w]_q$, where $w\in {\cal A}$, and $q$ is a primitive
$N$-th root of unity, is a ${\mathbb Z}_N$-graded $q$-differential
algebra and $d_v$ is
its $N$-differential.}%

Let us remind that a first order differential calculus over an
associative unital algebra $\cal B$ is a pair $({\cal M}, d)$,
where $\cal M$ is a $({\cal B},{\cal B})$-bimodule and $d$ is a
linear mapping $d:{\cal B}\to {\cal M}$ which satisfies the
Leibniz rule $d(ww')=d(w)w'+wd(w')$, where $\;w,w'\in {\cal B}$.
The subspace ${\cal A}^0$ of elements of grading zero of a
${\mathbb Z}_N$-graded algebra ${\cal A}$ is a subalgebra, and
$N$-differential $d_v$ restricted to this subalgebra induces a
first order differential calculus $({\cal A}^1, d_v)$ where the
space ${\cal A}^1$ of elements of grading 1 has a $({\cal
A}^0,{\cal A}^0)$-bimodule structure. Indeed it follows from the
associativity of the algebra $\cal A$ and its ${\mathbb
Z}_N$-graded structure that for each $k$ the mappings ${\cal
A}^0\times{\cal A}^k\to {\cal A}^k$ and ${\cal A}^k\times{\cal
A}^0\to{\cal A}^k$ determined by the algebra multiplication
$(r,w)\to rw,\;(w,s)\to ws$, where $w\in {\cal A}^k$ and $r,s\in
{\cal A}^0$, induce a $({\cal A}^0,{\cal A}^0)$-bimodule structure
on ${\cal A}^k$. In the next section we consider a reduced quantum
plane from a point of view of graded $q$-differential algebra and
study the first order differential calculus induced by the
$N$-differential.

\section{Reduced quantum plane as a $q$-differential algebra}
An exterior calculus with exterior differential $d$ satisfying
$d^N=0$ has been studied in
(\cite{Abramov1}.\cite{Abramov2},\cite{Abramov3}). In this section
we construct and study this kind of exterior calculus on a reduced
quantum plane with the help of the construction described in the
previous section. Let us remind that the unital associative
algebra ${\mathbb C}_{rq}$ generated, over the complex numbers
$\mathbb{C}$, by the two variables $x$ and $y$ satisfying the
relations
$x\,y=q\;y\,x,\;x^N=y^N={\bf 1},%
\label{relations of qp}$ where $q$ is a primitive $N$-th root of
unity and $\bf 1$ is the unity element of ${\mathbb C}_{rq}$, can
be considered as an algebra of polynomials over a reduced quantum
plane. Let us mention that this algebra has a representation  by
$N\times N$ complex matrices.

The set of monomials $B=\{{\bf
1},y,x,x^2,yx,y^2,\ldots,y^kx^l,\ldots,y^{N-1}x^{N-1}\}$ can be
taken as the basis of the vector space of the algebra
$\reducedqp$. Having chosen the basis $B$ we can endow this vector
space with a ${\mathbb Z}_N$-graded structure as follows: if a
polynomial $w\in \reducedqp$ written in terms of the monomials in
the basis $B$ has the form
\begin{equation}
w=\sum_{l=0}^{N-1} \beta_{l}y^kx^l,\quad \beta_l\in {\mathbb C}, %
        \quad k\in {\bf N},
\label{homogeneous polynomial}
\end{equation}
then we shall refer to it as the homogeneous polynomial of grading
${k}$, where ${k}\in {\mathbb Z}_N$. Let us denote the grading of
a homogeneous polynomial $w$ by $|w|$ and the subspace of the
homogeneous polynomials of grading $k$ by $\reducedqpk$. It is
obvious that
\begin{equation}
\reducedqp={\mathbb C}^{0}_{rq}\oplus {\mathbb C}^{ 1}_{rq}
  \oplus\ldots\oplus {\mathbb C}^{{N-1}}_{rq}.%
  \label{graded structure}
\end{equation}
In particular a polynomial $r$ of grading zero has the form
\begin{equation}
r=\sum_{l=0}^{N-1} \beta_lx^l,\qquad \beta_l\in {\mathbb
C},r\in{\mathbb C}^{0}_{rq}.
\end{equation}
It is easy to show that ${\mathbb Z}_N$-graded structure defined
by (\ref{graded structure})on a vector space $\reducedqp$ is
consistent with the algebra structure of $\reducedqp$, i.e. for
any two homogeneous polynomials we have $|ww'|=|w|+|w'|$.
Consequently $\reducedqp$ is a $Z_N$-graded algebra with respect
to (\ref{graded structure}), and there exists an element $v$ of
grading one of this algebra satisfying $v^N=\alpha\cdot {\bf 1}$,
where $\alpha\in {\mathbb C}$. Indeed one can take for instance
$v=y$ which satisfies all mentioned above conditions. According to
the first section we can endow a reduced quantum plane with the
structure of a graded $q$-differential algebra defining the
$N$-differential by the formula $d_vw=[v,w]_q$, where $q$ is a
primitive $N$-th root of unity and $w\in \reducedqp$.

As it was mentioned previously the subspace $\reducedqpzero$ of
elements of grading zero is a subalgebra of the algebra
$\reducedqp$. From a point of view of differential geometry we can
interpret the generator $x$ as a coordinate of a one-dimensional
space, the subalgebra $\reducedqpzero$ as an algebra of
(polynomial) functions or differential forms of degree 0 on this
one-dimensional space. Let us remind that the subalgebra
$\reducedqpzero$ is a commutative algebra generated by $x$
satisfying the single relation $x^N={\bf 1}$.

The subspace $\reducedqpk$ %
of polynomials of grading $k$ is a bimodule over the algebra of
functions $\reducedqpzero$. If we put a homogeneous polynomial $w$
of grading $k$ to the form
\begin{equation}
w=y^k\sum_{l=0}^{N-1} \beta_lx^l=y^k\,r,\quad%
r=\sum_{l=0}^{N-1}\beta_lx^l\in\reducedqpzero,
\end{equation}
and take into account that the polynomial $r=(y^k)^{-1}w=y^{N-k}w$
is uniquely determined then we can conclude that $\reducedqpk$ is
a free right module over $\reducedqpzero$ generated by $y^k$. Thus
extending our differential-geometric interpretation to the whole
algebra $\reducedqp$ we can interpret the bimodule $\reducedqpk$
as a module of differential forms of degree $k$ over the algebra
of functions $\reducedqpzero$ and $d_v$ as an exterior
differential.

It is well known that a bimodule structure on a free right module
is uniquely determined by the homomorphism of the corresponding
algebra, and in our case this means that there exists a
homomorphism $A_{k}:\reducedqpzero\to\reducedqpzero$ such that
\begin{equation}
r\,y^k=y^k\,A_{k}(r),\qquad r\in\reducedqpzero.%
\label{bimodule srtucture}
\end{equation}
It is easy to find that $A_{k}(x)=q^k\,x$ and $A_k=A_1^k,A_0=I$
where $I:\reducedqpzero\to \reducedqpzero$ is the identity
mapping. Thus we have the set $\{A_{k}\}_{k=0}^{N-1}$ of
homomorphisms of the algebra $\reducedqpzero$.

Since $\reducedqpone$ is a free right module over $\reducedqpzero$
there exists an invertible element $u\in \reducedqpzero$ such that
$v=y\cdot u$. We can take $v$ for a generator of the free right
module $\reducedqpone$. Using the relation (\ref{bimodule
srtucture}) and the commutativity of the algebra $\reducedqpzero$
we find that the relation determining the bimodule structure in
terms of the new generator $v$ has the same form (\ref{bimodule
srtucture}) as in the case of the generator $y$. Now we can write
the polynomial $d_vw$ in the form of an element of the right
module generated by $v$ as follows
\begin{eqnarray}
d_vw=vw-wv=vw-vA_1(w)=v(w-A_1(w))=v\,\Delta_0(w).%
\end{eqnarray}
where $\Delta_0=I-A_1:\reducedqpzero\to \reducedqpzero$. It is
easy to check that given any polynomials $w,w'\in \reducedqpzero$
the mapping $\Delta_0$ satisfies the following properties
\begin{eqnarray}
\Delta_0(ww')&=&\Delta_0(w)w'+A_1(w)\Delta_0(w'),%
                           \label{property1 of delta}\\%
\Delta_0(x^k)&=&(1-q)[k]_q\;x^k.%
    \label{propert2 of delta}
\end{eqnarray}
This formula shows that $d_vx$ can be taken as a generator of the
right module $\reducedqpone$.

It is well known \cite{Borowiec_Kharchenko1} that a differential
on a unital associative algebra induces the right partial
derivatives which satisfy a generalized Leibniz rule. In our case
we have only one derivative
$\partial:\reducedqpzero\to\reducedqpzero$ which is defined by the
formula $d_vw=d_vx\,\,\partial w,\; \forall w\in \reducedqpzero$.
The explicit formula for the partial derivative has the form
$\partial w=(1-q)^{-1}x^{N-1}\Delta_0(w)$. It follows immediatly
from (\ref{property1 of delta}),(\ref{propert2 of delta}) that
this derivative satisfies the generalized Leibniz rule
$\partial(ww')=\partial(w)\cdot w'+A_1(w)\cdot\partial(w')$ and %
$\partial x^k=[k]_q\;x^{k-1}$.

The author thanks the organizers of the "International Conference
on High Energy and Mathematcal Physics" for their hospitality and
acknowledges the financial support by grant ETF 6206 of the
Estonian Science Foundation.

\bibliographystyle{amsplain}

\end{document}